\documentclass[reqno]{amsart}
\usepackage{amsmath,amssymb,graphicx,amscd}
\usepackage{mathrsfs}
\usepackage{mathrsfs}

\DeclareMathOperator{\E}{{\mathbb{E}}}

\newtheorem{thm}{Theorem}[section]

\newtheorem{lem}[thm]{Lemma}
\newtheorem{prop}[thm]{Proposition}
\theoremstyle{definition}
\newtheorem{defn}[thm]{Definition}
\theoremstyle{remark}
\newtheorem{rem}[thm]{Remark}
\theoremstyle{remark}

\numberwithin{equation}{section}

\begin{document}

\title{A generalization of Doob's maximal identity}

\author{Ashkan Nikeghbali}
 \address{Institut f\"ur Mathematik, Universit\"at Z\"urich, Winterthurerstrasse 190, CH-8057 Z\"urich, Switzerland}
 \email{ashkan.nikeghbali@math.unizh.ch}

 \subjclass[2000]{05C38, 15A15, 15A18} \keywords{Progressive enlargements of filtrations,
Az\'{e}ma's supermartingale, Maximum of local martingales
General theory of stochastic processes,  Last
passage time, Multiplicative decompositions}

\date{\today}

\begin{abstract}
In this paper, using martingale techniques, we prove a generalization of Doob's maximal identity in the setting of  continuous nonnegative local submartingales  $\left(X_{t}\right)$ of the form:
$X_{t}=N_{t}+A_{t}$, where the measure $\left(dA_{t}\right)$ is
carried by the set $\left\{t:\;X_{t}=0\right\}$.  In particular, we give a multiplicative decomposition for the Az\'ema supermartingale associated with some last passage times related to such processes and we prove that these non-stopping times contain very useful information.   As a consequence, we obtain the law of the maximum of a continuous nonnegative local martingale $(M_t)$ which satisfies $M_\infty=\psi(\sup_{t\geq0}M_t)$ for some measurable function $\psi$ as well as the law of the last time this maximum is reached. 

\end{abstract}

\maketitle

\section{Introduction}

Throughout this paper, $\left(\Omega,\mathcal{F},\left(\mathcal{F}_{t}\right),\mathbb{P}\right)$ will 
be a filtered probability space, satisfying the usual assumptions. If $(X_t)$ denotes a stochastic process, we note $\overline{X}_t=\sup_{s\leq t}X_s$. We also recall the definition of local submartingales of the class $(\Sigma)$ which will play a crucial role in this paper. These processes were first introduced by Yor (\cite{yorinegal}) and further studied in \cite{AshkanYorremII,AshkanYorrem}. For simplicity, we restrict ourselves to the case of continuous processes:
\begin{defn}\label{martreflechies}
Let $\left(X_{t}\right)$ be a positive local submartingale, which
decomposes as:
$$X_{t}=N_{t}+A_{t}.$$
We say that $\left(X_{t}\right)$ is of class $(\Sigma)$ if:
\begin{enumerate}
\item $\left(N_{t}\right)$ is a continuous local martingale, with $N_{0}=0$;
\item $\left(A_{t}\right)$ is a continuous increasing process, with $A_{0}=0$;
\item the measure $\left(dA_{t}\right)$ is carried by the set
$\left\{t:\;X_{t}=0\right\}$.
\end{enumerate}If additionally, $\left(X_{t}\right)$
is of class $(D)$, we shall say that $\left(X_{t}\right)$ is of
class $(\Sigma D)$.
\end{defn}
Well known examples of such
stochastic processes are $\left(\left|M_{t}\right|\right)$, the
absolute value of some local martingale $M$, $M^+$ or
$\left(\overline{M}_{t}-M_{t}\right)$. Many other processes or simple transforms of diffusions fall in the class $(\Sigma)$ (see \cite{AshkanYorrem} for more details). 

\begin{rem}
In \cite{AshkanYorrem}, the local martingale part $N_t$ is allowed to be right continuous with left limit. Here, for simplicity, we restrict our attention to the case of continuous processes. But some of the results in this paper would hold in the more general setting under some extra assumptions (such as conditions on  the sign of the jumps). 
\end{rem}

In \cite{ashyordoob}, nonnegative local martingales $(M_t)$ which satisfy $\lim_{t\to\infty}M_t=0$ are studied in depth: in particular, the law of $\sup_{t\geq0}M_t$ (which is distributed as $\dfrac{M_0}{U}$, where $U$ is a uniform random variable, and known as Doob's maximal identity\footnote{In fact, this identity, which is a simple consequence of the optional stopping theorem, has been  well known for some time and appears for example as an exercise in \cite{revuzyor} p.73.}) as well as the random time $g=\sup\{t: M_t=\overline{M}_t\}$,  are shown to play an important role in the characterization of honest times (i.e. ends of predictable sets) and in the theory of progressive enlargements of filtrations (see also \cite{columbia}). One key formula in \cite{ashyordoob} is the multiplicative decomposition of the Az\'ema's supermartingale associated with $g$:
\begin{equation}\label{eq:g}
\mathbb{P}[g>t|\mathcal{F}_{t}]=\dfrac{M_t}{\overline{M}_t}.
\end{equation}
One may wonder if a simple multiplicative decomposition such as (\ref{eq:g}) still holds in more general situations?  One such situation occurs in the resolution of the Skorokhod stopping problem by Az\'ema and Yor  (\cite{AY}, see for example the survey \cite{jansurvey} for an overview) where $M_\infty=\psi(\overline{M}_\infty)$, with $\psi$  some measurable function. 

The goal of this paper is twofold:
\begin{itemize}
\item to show that the multiplicative decomposition formula  (\ref{eq:g}) as well Doob's maximal identity naturally extend in the framework of local submartingales of the class $(\Sigma)$;
\item to outline the importance of honest times in the theory of stochastic processes: it is a remarkable fact that the knowledge of the Az\'ema supermartingales (as well as the law) of such random times gives  information about the underlying process.
\end{itemize}
More precisely, this paper is organized as follows:

In Section 2.1, we recall some useful results about honest times  as well as an elementary additive characterization of Az\'ema's supermartingale associated with the end of a predictable set. 

In Section 2.2, we consider conveniently stopped processes $(X_{t\wedge T})$ of the class $(\Sigma)$, with possibly $T=\infty$, and which satisfy  $X_T=\psi(A_T)$, for some measurable function $\psi$. We show how the results of section 2.1 can be used to obtain a multiplicative decomposition for the Az\'ema's supermartingale associated with  the last zero of $(X_{t\wedge T})$ and use this decomposition to compute the law of $A_T$. As a corollary, we obtain a generalization of the results in \cite{ashyordoob} for nonnegative local martingales which satisfy $M_\infty=\psi(\overline{M}_\infty)$.

\section*{Acknowledgements}
I am very grateful to Marc Yor for all his comments and for giving me access to  preliminary versions of \cite{mry}. The present paper is extracted from my PhD thesis \cite{phdashkan}, page 253-260.

\section{A generalization of Doob's maximal identity}

\subsection{Useful results about honest times}
Let
$\left(\Omega,\mathcal{F},\left(\mathcal{F}_{t}\right),\mathbb{P}\right)$
be a filtered probability space, satisfying the usual assumptions,
and let $L$ be the end of an
$\left( \mathcal{F}_{t}\right) $ predictable set $%
\Gamma $, i.e:
\begin{equation*}
L=\sup \left\{ t:\left( t,\omega \right) \in \Gamma \right\}.
\end{equation*} A typical example of such a random time $L$ is $$L=\sup\left\{t\leq1:\; B_{t}=0\right\},$$ where as usual
$\left(B_{t}\right)$ denotes the standard Brownian Motion. Ends of predictable sets are the most studied random times after stopping times (see \cite{Ashkanessay} for examples and more references).
One process plays an essential role in
the study of the random time $L$, namely the supermartingale:
\begin{equation}
    Z_{t}^{L}=\mathbb{P}\left(L>t|\mathcal{F}_{t}\right),
\end{equation}
associated with $L$ by Az\'ema in \cite{azema}, and chosen to be
c\`{a}dl\`{a}g. 
For simplicity, we make the following  \underline{assumptions throughout this
paper, which we call the} $\mathbf{(CA)}$ conditions:
\begin{enumerate}
\item all $\left( \mathcal{F}_{t}\right) $-martingales are \textbf{\underline{c}}ontinuous (e.g:
 the Brownian filtration).

\item the random time $L$ \textbf{\underline{a}}voids every $\left( \mathcal{F}_{t}\right) $%
-stopping time $T$, i.e. \\ 
$\mathbb{P}\left[ L =T\right] =0$.\bigskip
\end{enumerate}
\begin{rem}
Under the $\mathbf{(CA)}$ conditions, the optional and the
predictable sigma fields are equal and the supermartingale
$(Z_{t}^{L})$ is continuous.
\end{rem}

Now consider the Doob-Meyer
decomposition of $Z^{L}$:
\begin{equation}\label{doobmeyer}
Z_{t}^{L}=1+\mu _{t}^{L}-A_{t}^{L}
\end{equation}
The process $\left( A_{t}^{L}\right)$, which we shall sometimes note
$\left( A_{t}\right)$ in the sequel, is the dual predictable
projection of the increasing process $\mathbf{1}_{\left\{ L\leq
t\right\} }$, and
\begin{equation*}
\mu_{t}^{L}=\mathbb{E}\left( A_{\infty }^{L}\mid
\mathcal{F}_{t}\right) -1
\end{equation*}
\begin{prop}[Az\'{e}ma \cite{azema}]\label{lemazema}
Let $L$ be the end of some predictable set; then
$$L=\sup\left\{t:\;Z_{t}=1\right\},$$and the measure $dA_{t}$ is
carried by the set $\left\{t:\;Z_{t}=1\right\}$. In particular, $A$
does not increase after $L$, i.e. $A_{L}=A_{\infty}$. Moreover,  $A_{\infty }$ is
distributed with the standard exponential law.
\end{prop}
It follows from this proposition that $$X_{t}\equiv 1-Z_{t}^{L}$$ is
of the class $(\Sigma D)$, with: $X_{0}=0$ and
$\displaystyle\lim_{t\rightarrow\infty}X_{t}=1$. In fact, the converse of this result is also true. It is stated in \cite{AshkanYorremII} (where it was used to derive path decomposition results for diffusion processes) as a consequence of  a very useful lemma (although not so well known) which first appeared  in the papers of Az\'ema, Meyer and Yor
\cite{azemameyeryor} (section 7) and Az\'ema and Yor \cite{azemayorzero} (Proposition 2.2), in a
more abstract framework. It will again play an important role in our discussions in this paper. We  state a  variant of this lemma  in which the local martingale $(N_t)$   in Definition \ref{martreflechies} does not necessarily start from $0$ but from an arbitrary positive number $x$: $N_0=x\geq0$. The proof by Az\'ema and Yor \cite{azemayorzero} (Proposition 2.2) still holds  and we   simply reproduce it.

\begin{lem}\label{martingalerefelchies}
Let $\left(X_{t}\right)$ be a submartingale of the class $(\Sigma
D)$ with $X_0=x\geq0$ (i.e. $N_0=x$ in Definition \ref{martreflechies}) and let
$$L=\sup\left\{t:\;X_{t}=0\right\},$$ with the convention that $ \sup\{\emptyset\}=0$. Assume further that:
$$\mathbb{P}\left(X_{\infty}=0\right)=0.$$ Then:
\begin{equation}\label{martrel}
    X_{t}=\mathbb{E}\left(X_{\infty}\mathbf{1}_{\left\{L\leq
    t\right\}}|\mathcal{F}_{t}\right).
\end{equation}

\end{lem}
\begin{proof}
Since $(X_{t})$ is continuous, the set $\left\{t:\;X_{t}=0\right\}$
is a predictable closed set. Let us remark that:
$$X_{\infty}\mathbf{1}_{\left\{L\leq t\right\}}=X_{d_{t}},$$where
$$d_{t}\equiv \inf\left\{s>t:\;X_{s}=0\right\}.$$ Hence, we have:
\begin{eqnarray*}
  \mathbb{E}\left(X_{\infty}\mathbf{1}_{\left\{L\leq t\right\}}|\mathcal{F}_{t}\right) &=& \mathbb{E}\left(X_{d_{t}}|\mathcal{F}_{t}\right)\\
   &=&
   \mathbb{E}\left(N_{d_{t}}|\mathcal{F}_{t}\right)+\mathbb{E}\left(A_{d_{t}}|\mathcal{F}_{t}\right).
\end{eqnarray*} Now, from the optional stopping theorem, we have:
$$\mathbb{E}\left(N_{d_{t}}|\mathcal{F}_{t}\right)=N_{t}.$$
Moreover, as $\left(dA_{t}\right)$ is carried by the set
$\left\{t:\;X_{t}=0\right\}$, we have:$$A_{d_{t}}=A_{t}.$$ We can
thus conclude
that:$$\mathbb{E}\left(X_{\infty}\mathbf{1}_{\left\{L\leq
t\right\}}|\mathcal{F}_{t}\right)=N_{t}+A_{t}=X_{t},$$and this
completes the proof.
\end{proof}

\begin{rem}
In their  recent papers, D. Madan, B. Roynette and M. Yor (\cite{mry,my})  have obtained some nice and striking results relating the price of some options to the Az\'ema supermartingale associated with  some last passage times.  The very first step of this work can  be viewed through
 Lemma \ref{martingalerefelchies}. Indeed, consider $(M_t)$ a continuous and  nonnegative local martingale starting from $z$ and such that $\lim_{\to\infty}M_t=0$. In the pricing of put options, quantities such as $\E[(K-M_T)^+]$ or  $\E[(K-M_T)^+|\mathcal{F}_t]$ for some $t\leq T$ often arise.
 
 For any $K>0$, consider the local submartingale $X_t=(K-M_t)^+$. It is easy to see that $X$ satisfies the assumptions of Lemma \ref{martingalerefelchies} with $x=(K-z)^+$ and $X_\infty=K$. Here, $L\equiv g_K=\sup\{t:\;M_t=K\}$. An application of Lemma \ref{martingalerefelchies}  yields (these expressions first appeared in \cite{mry,my}):
 $$(K-M_t)^+=K\mathbb{P}[g_K\leq t|\mathcal{F}_t].$$Consequently, for $t\leq T$, we have:
$$\E[(K-M_T)^+|\mathcal{F}_t]=K\mathbb{P}[g_K\leq T|\mathcal{F}_t],$$and $$\E[(K-M_T)^+]=K\mathbb{P}[g_K\leq T].$$
For other proofs and more results with financial models in view, the interested reader should refer to \cite{mry,my}. 
 
Moreover, at $t=0$, the above representation  gives:
$$\mathbb{P}[g_K=0]=\dfrac{(K-z)^+}{K}.$$ But $\{g_K=0\}=\{\overline{M}_\infty<K\}$ and hence $$\mathbb{P}[\overline{M}_\infty<K]=\dfrac{(K-z)^+}{K},$$which is Doob's maximal identity. More generally, if $\left(X_{t}\right)$ satisfies the assumptions of Lemma \ref{martingalerefelchies} with $X_\infty=K$ where $K>0$ is a constant, then $\mathbb{P}[\overline{X}_\infty>0]=\dfrac{x}{K}$.
\end{rem}
\begin{rem}
In a forthcoming paper \cite{NPY}, some other applications of honest times in financial modelling are developed.
\end{rem}
\begin{rem}
More generally, as already mentioned, Lemma \ref{martingalerefelchies} leads to a very natural characterization of Az\'ema's supermartingales: roughly speaking, any submartingale $X$ of the class ($\Sigma$) that converges to a nonzero constant $K$, that we can without loss of generality take to be one, is equal to $1-\mathbb{P}[L>t|\mathcal{F}_t]$, where $L$ is the last time $X$ hits zero. This characterization corresponds to Theorem 2.3 in \cite{AshkanYorremII} and is systematically used there to express some well known submartingales of the class $(\Sigma)$ in terms of some associated last passage times.
\end{rem}
\subsection{Az\'ema's supermartingale associated with the last zero}

In \cite{AshkanYorrem}, we computed the law of $A_{\infty}$ when
$X$ is of the class $(\Sigma D)$; in particular, we were interested
in stopping times of the form
\begin{equation}\label{tpasarret}
T\equiv\inf\left\{t:\;\varphi\left(A_{t}\right)X_{t}\geq1\right\},
\end{equation}for certain nonnegative Borel functions.  In particular, we
observed in \cite{AshkanYorrem} that
$\left(\varphi\left(A_{t}\right)X_{t}\right)$ is of the class
$(\Sigma)$, and may be represented as:
\begin{equation}\label{represdepitimex}
\varphi\left(A_{t}\right)X_{t}=\int_{0}^{t}\varphi\left(A_{u}\right)dN_{u}+\Phi\left(A_{t}\right),
\end{equation}where $\Phi\left(x\right)=\int_{0}^{x}dz\varphi\left(z\right)$ (since $X$ is continuous, this last formula can also be easily
deduced from standard balayage arguments as exposed in
\cite{revuzyor}, chapter VI). We also showed that if $\varphi$ is a
nonnegative locally bounded Borel function, such that
$\int_{0}^{\infty}dz\varphi\left(z\right)=\infty$, then the stopping
time $T$ is almost surely finite:
\begin{equation}\label{finitudedeT}
\int_{0}^{\infty}dz\varphi\left(z\right)=\infty\Rightarrow T<\infty.
\end{equation}

Here, we shall develop further the study  of  $X$ on $[0,T]$, using
Lemma \ref{martrel} and Proposition \ref{lemazema}. More precisely,
we shall look for the  Az\'{e}ma's supermartingale  associated with the random time
$$g_{T}\equiv\sup\left\{t<T:\;X_{t}=0\right\}.$$
\begin{prop}\label{loicondegt}
Let $\varphi$ be a nonnegative locally bounded  Borel function, such
that $\varphi\left(x\right)>0,\;\forall x>0$ and such that
$\int_{0}^{\infty}dz\varphi\left(z\right)=\infty$ and let $X$ be a
local submartingale of the class $(\Sigma)$, such that
$A_{\infty}=\infty$. Let $T$ be defined as in (\ref{tpasarret}) and
let
$$g_{T}\equiv\sup\left\{t<T:\;X_{t}=0\right\}.$$Then,
\begin{equation}\label{surmartdegt}
\mathbb{P}\left(g_{T}\leq
t|\mathcal{F}_{t}\right)=\varphi\left(A_{t\wedge T}\right)X_{t\wedge
T}=\int_{0}^{t\wedge
T}\varphi\left(A_{u}\right)dN_{u}+\Phi\left(A_{t\wedge T}\right),
\end{equation}and $g_{T}$ avoids all $\left(\mathcal{F}_{t}\right)$
stopping times. Consequently,
$$\mathbb{P}\left(A_{T}>x\right)=\exp\left(-\int_{0}^{x}dz\varphi\left(z\right)\right).$$
\end{prop}
\begin{proof}
From (\ref{finitudedeT}), $T<\infty$, and consequently,
$\varphi\left(A_{T}\right)X_{T}=1$. Now, since
$\left(\varphi\left(A_{t\wedge T}\right)X_{t\wedge T}\right)$ is of
the class $(\Sigma)$, we can apply Lemma \ref{martrel} (with
$X_{\infty}=1$) to obtain the identity (\ref{surmartdegt}). From
(\ref{represdepitimex}), we also have: $\mathbb{P}\left(g_{T}\leq
t|\mathcal{F}_{t}\right)=\int_{0}^{t\wedge
T}\varphi\left(A_{u}\right)dN_{u}+\Phi\left(A_{t\wedge T}\right)$,
and hence the dual predictable projection of
$\left(\mathbf{1}_{g_{T}\leq t}\right)$ is
$\left(\Phi\left(A_{t\wedge T}\right)\right)$. Now, since all
$\left(\mathcal{F}_{t}\right)$ martingales are continuous and
$\left(\Phi\left(A_{t\wedge T}\right)\right)$ is continuous, all
$\left(\mathcal{F}_{t}\right)$ stopping times are predictable, and
$g_{T}$ avoids all $\left(\mathcal{F}_{t}\right)$
stopping times. Indeed, for any $\left( \mathcal{F}_{t}\right) $
stopping time $R$,
\begin{equation*}
\mathbf{E}\left[ \mathbf{1}_{\left\{ g_{T}=R\right\}
}\right]=\mathbb{E}\left[\Delta \left(\Phi\left(A_{R\wedge
T}\right)\right)\right]=0.
\end{equation*}
Thus we get $P\left( g_{T}=R\right) =0$.

The fact that
$\mathbb{P}\left(A_{T}>x\right)=\exp\left(-\int_{0}^{x}dz\varphi\left(z\right)\right)$
is a consequence of the fact that $\Phi\left(A_{T}\right)$ is
distributed as a random variable with the standard exponential law
(see Proposition \ref{lemazema}).
\end{proof}
\begin{rem}
Applying Corollary 4.2 in \cite{ashyordoob}, one obtains:
$\mathbb{P}\left(g_{T}>
t|\mathcal{F}_{t}\right)=\dfrac{Y_t}{\overline{Y}_t}$, with  $(Y_t)$ some continuous and nonnegative local martingale starting from $1$ and converging to $0$ at infinity. Moreover, $Y_t$ and $\overline{Y}_{t}$ can be expressed in terms of $N_t$ and $A_t$:
\begin{eqnarray*}
  Y_{t} &=& \exp \left( \int_{0}^{t}\frac{dM_{s}}{Z_{s}}-\frac{1}{2}\int_{0}^{t}%
\frac{d<M>_{s}}{Z_{s} ^{2}}\right) \\
  \overline{Y}_{t} &=& \exp\left(\Phi\left(A_{t\wedge T}\right)\right);
\end{eqnarray*}with 
$M_t=1+\int_{0}^{t\wedge
T}\varphi\left(A_{u}\right)dN_{u}$ and $Z_{t}=\varphi\left(A_{t\wedge T}\right)X_{t\wedge
T}$. The multiplicative expression (\ref{surmartdegt}) is thus much simpler.
\end{rem}
It may happen that $T=\infty$, or in other words,
$X_{\infty}=\varphi\left(A_{\infty}\right)$, for some nonnegative
Borel functions, in which case we have the following proposition:
\begin{prop}\label{loicinddeginf}
Let $\varphi$ be a nonnegative Borel function such that $1/\varphi$
is locally bounded, and let $X$ be a submartingale of the class
$(\Sigma D)$, such that $X_{\infty}=\varphi\left(A_{\infty}\right)$.
Define: $$g\equiv\sup\left\{t:\;X_{t}=0\right\}.$$If
$\varphi\left(A_{\infty}\right)>0,\;a.s.$, then:
\begin{equation}\label{surmartdeginfi}
    \mathbb{P}\left(g\leq
t|\mathcal{F}_{t}\right)=\dfrac{X_{t}}{\varphi\left(A_{t}\right)}=\int_{0}^{t}\dfrac{dN_{u}}{\varphi\left(A_{u}\right)}+\int_{0}^{A_{t}}\dfrac{dz}{\varphi\left(z\right)}.
\end{equation}Consequently,
$$\mathbb{P}\left(A_{\infty}>x\right)=\exp\left(-\int_{0}^{x}\dfrac{dz}{\varphi\left(z\right)}\right).$$
\end{prop}
\begin{proof}
From Lemma \ref{martrel}, we have:
$$X_{t}=\mathbb{E}\left(X_{\infty}\mathbf{1}_{\left\{g\leq
    t\right\}}|\mathcal{F}_{t}\right)=\mathbb{E}\left(\varphi\left(A_{\infty}\right)\mathbf{1}_{\left\{g\leq
    t\right\}}|\mathcal{F}_{t}\right).$$But from Proposition
    \ref{lemazema}, we also have:
    $$\varphi\left(A_{\infty}\right)\mathbf{1}_{\left\{g\leq
    t\right\}}=\varphi\left(A_{t}\right)\mathbf{1}_{\left\{g\leq
    t\right\}},$$and consequently,
    $$X_{t}=\varphi\left(A_{t}\right)\mathbb{P}\left(g\leq
t|\mathcal{F}_{t}\right),$$ and (\ref{surmartdeginfi}) follows
easily from (\ref{represdepitimex}). Eventually, the law of
$A_{\infty}$ follows from the fact that the dual predictable
projection of $\mathbf{1}_{\left\{g\leq
    t\right\}}$, which is $\int_{0}^{A_{t}}\frac{dz}{\varphi\left(z\right)}$, taken at $t=\infty$, follows the standard exponential law.
\end{proof}
\begin{rem}
We note that taking $\varphi\equiv K>0$ a constant, we recover the results discussed in the remarks following Lemma \ref{martingalerefelchies}. We moreover see that $A_\infty$ follows an exponential law with parameter $1/K$.
\end{rem}

\bigskip

Now, we give some applications of Propositions \ref{loicondegt} and
\ref{loicinddeginf} to the case when $M_{\infty}=\psi\left(\overline{M}_{\infty}\right)$, with $\psi$ a
nonnegative Borel function. The next Theorem generalizes Lemma 2.1 and Proposition 2.2 in \cite{ashyordoob}. Without loss of generality, we can assume that
$M_{0}=1$. With the notation
$$\overline{M}_{t}\equiv\sup_{u\leq t}M_{u},$$we have:
\begin{thm}\label{loimaxmlocpositifs}
Let $\left(M_{t}\right)$ be a continuous nonnegative local
martingale starting from $1$, and such that
$M_{\infty}=\psi\left(\overline{M}_{\infty}\right)$, with $\psi$ a
nonnegative Borel function such that
$\psi\left(x\right)<x,\;\forall x\geq1$. 
Define
$$g\equiv\sup\left\{t:\;M_{t}=\overline{M}_{t}\right\}.$$
Then,
\begin{eqnarray*}
 \mathbb{P}\left(g\leq
t|\mathcal{F}_{t}\right) &=& \dfrac{\overline{M}_{t}-M_{t}}{\overline{M}_{t}-\psi\left(\overline{M}_{t}\right)} \\
   &=&
   -\int_{0}^{t}\dfrac{dM_{u}}{\overline{M}_{u}-\psi\left(\overline{M}_{u}\right)}+\int_{1}^{\overline{M}_{t}}\dfrac{dz}{z-\psi\left(z\right)}.
\end{eqnarray*}Hence, the dual predictable
projection of the raw increasing process
$\left(\mathbf{1}_{\left\{g\leq
    t\right\}}\right)$, is $\left(\int_{0}^{\overline{M}_{t}}\frac{dz}{z-\psi\left(z\right)}\right)$. Consequently, we have:
\begin{equation}\label{generalisationdeidentitedoob}
    \mathbb{P}\left(\overline{M}_{\infty}>x\right)=\exp\left(-\int_{1}^{x}\dfrac{dz}{z-\psi\left(z\right)}\right),\;\forall x\geq1.
\end{equation}In particular, when $\psi\equiv0$, i.e. when $M_{\infty}=0$, we recover the
following results  which appear in
\cite{ashyordoob}:\begin{eqnarray*}
  \mathbb{P}\left(g>
t|\mathcal{F}_{t}\right) &=& \dfrac{M_{t}}{\overline{M}_{t}}, \\
  \mathbb{P}\left(\overline{M}_{\infty}>x\right)= &=&
  \dfrac{1}{x},\;\forall x\geq1.
\end{eqnarray*}
\end{thm}
\begin{proof}
Let us define
$$X_{t}\equiv 1-\dfrac{M_{t}}{\overline{M}_{t}}.$$An application of
It\^{o}'s formula, combined with the fact that $d\overline{M}_{t}$
is carried by the set $\left\{t:\;M_{t}=\overline{M}_{t}\right\}$,
yields:
$$X_{t}=-\int_{0}^{t}\dfrac{dM_{u}}{\overline{M}_{u}}+\log\left(\overline{M}_{t}\right),$$and
since
$\left\{t:\;X_{t}=0\right\}=\left\{t:\;M_{t}=\overline{M}_{t}\right\}$,
we easily deduce that $X$ (which is bounded by $1$) is of the class
$(\Sigma D)$. With the notations of Proposition
\ref{surmartdeginfi}, we have
$A_{t}=\log\left(\overline{M}_{t}\right)$ and
$\varphi\left(x\right)\equiv1-\dfrac{\psi\left(\exp\left(x\right)\right)}{\exp\left(x\right)}$.
Since $\psi\left(x\right)<x$, we have
$\varphi\left(A_{\infty}\right)>0$, and the results of the corollary
follow from an application of Proposition \ref{loicinddeginf} and
some elementary calculations.
\end{proof}
\begin{rem}
The situation described in the previous Theorem often occurs in
the Az\'{e}ma-Yor
solution to Skorokhod's embedding problem which relies upon the
construction of a Brownian martingale $M_{t}=B_{t\wedge T}$ such
that $M_{T}=\psi(\overline{M}_{T})$ (see \cite{AY}).
\end{rem}


\begin{thebibliography}{99}

\bibitem{azema} \textsc{J. Az\'{e}ma}: \textit{Quelques applications de la th%
\'{e}orie g\'{e}n\'{e}rale des processus I}, Invent. Math.
\textbf{18} (1972) 293-336.


\bibitem{azemameyeryor} \textsc{J. Az\'ema, P.A. Meyer, M. Yor}: \textit{Martingales relatives}, S\'{e}m.Proba. XXVI, Lecture Notes in Mathematics \textbf{1526},
(1992), 307-321.

\bibitem{AY} \textsc{J. Az\'{e}ma, M. Yor}: \textit{Une solution simple au
probl\`{e}me de Skorokhod}, S\'{e}m.Proba. XIII, Lecture Notes in
Mathematics \textbf{721}, (1979), 90-115 and 625-633.


\bibitem{azemayorzero} \textsc{J. Az\'ema, M. Yor}: \textit{Sur les z\'{e}ros des martingales continues}, S\'{e}m.Proba. XXVI, Lecture Notes in Mathematics \textbf{1526},
(1992), 248-306.

\bibitem{barlow} \textsc{M.T. Barlow}, \textit{Study of a filtration
expanded to include an honest time}, ZW, \textbf{44}, 1978, 307-324.


\bibitem{dellachmeyer} \textsc{C. Dellacherie, P.A. Meyer}: \textit{%
Probabilit\'{e}s et potentiel}, Hermann, Paris, volI. 1980.

\bibitem{delmaismey} \textsc{C. Dellacherie, B. Maisonneuve, P.A. Meyer}:
\textit{Probabilit\'{e}s et potentiel}, Chapitres XVII-XXIV:
Processus de Markov (fin), Compl\'{e}ments de calcul stochastique,
Hermann (1992).


\bibitem{jeulin} \textsc{T. Jeulin}: \textit{Semi-martingales et
grossissements d'une filtration}, Lecture Notes in Mathematics
\textbf{833}, Springer (1980).


\bibitem{jeulinyor81} \textsc{T. Jeulin, M. Yor}: \textit{Sur les distributions de certaines fonctionnelles du mouvement brownien}, S\'{e}m.Proba. XV, Lecture Notes in Mathematics \textbf{850},
(1981), 210-226.

\bibitem{jeulinyor} \textsc{T. Jeulin, M. Yor (eds)}: \textit{Grossissements
de filtrations: exemples et applications}, Lecture Notes in
Mathematics \textbf{1118}, Springer (1985).

\bibitem{mry}  \textsc{D. Madan, B. Roynette, M. Yor}: \textit{From Black-Scholes formula, to local
times and last passage times for certain
submartingales}, preprint (September 2007).

\bibitem{my}  \textsc{D. Madan, B. Roynette, M. Yor}  \textit{Option prices as probabilities}, submitted (January 2008).

\bibitem{columbia} \textsc{R. Mansuy, M. Yor}: \textit{Random times and (enlargement of filtrations) in a Brownian
setting}, Lecture Notes in Mathematics, \textbf{1873}, Springer
(2006).

\bibitem{phdashkan} \textsc{A. Nikeghbali}: \textit{Temps al\'eatoires, filtrations et sousmartingales: quelques d\'eveloppements r\'ecents}, th\`ese de doctorat de l'universit\'e Paris 6, September 2005.

\bibitem{AshkanYorremII} \textsc{A. Nikeghbali}: \textit{Enlargements of filtrations and path decompositions at non stopping times}, Prob. Theory Related
Fields, \textbf{136} (4), 2006, 524-540.

\bibitem{ashyordoob} \textsc{A. Nikeghbali, M. Yor}: \textit{Doob's maximal identity, multiplicative decompositions
 and enlargements of filtrations}, Illinois Journal of Mathematics, \textbf{50}
 (4) 791-814 (2006).

\bibitem{AshkanYorrem} \textsc{A. Nikeghbali}: \textit{A class of remarkable submartingales}, Stochastic Processes and their Applications \textbf{116}, (2006) 917-938.



\bibitem{Ashkanessay} \textsc{A. Nikeghbali}: \textit{An essay on the general theory of stochastic processes}, Probab. Surv., \textbf{3}, 2006, 345-412.



\bibitem{NPY} \textsc{A. Nikeghbali, E. Platen}: \textit{On honest times in financial modelling}, in preparation.

\bibitem{jansurvey} \textsc{J. Ob\l\'{o}j}: \textit{The Skorokhod embedding
problem and its offspring}, Probab. Surv., \textbf{1}, (2004),
321-390.

\bibitem{pitmanyor} \textsc{J.W. Pitman, M. Yor}: \textit{Bessel processes
and infinitely divisible laws, In: D. Williams (ed.) Stochastic
integrals}, Lecture Notes in Mathematics \textbf{851}, Springer
(1981).


\bibitem{protter} \textsc{P.E. Protter}: \textit{Stochastic integration and
differential equations}, Springer. Second edition (2005).

\bibitem{revuzyor} \textsc{D. Revuz, M. Yor}: \textit{Continuous martingales
and Brownian motion}, Springer. Third edition (1999).

\bibitem{rogers} \textsc{L.C.G. Rogers}, \textit{The joint law of the maximum and terminal value of a martingale},
Prob. Theory Related Fields, \textbf{95}(4), (1993), 451-466.

\bibitem{rogerswilliams} \textsc{C. Rogers, D. Williams}: \textit{%
Diffusions, Markov processes and Martingales, vol 2: Ito calculus},
Wiley and Sons, New York, (1987).

\bibitem{vallois} \textsc{P. Vallois}: \textit{Sur la loi du maximum et du temps local d'une martingale continue uniform\'{e}ment int\'{e}grable},
Proc. London Math. Soc. \textbf{3}, 69(2) (1994), 399-427.


\bibitem{yorinegal} \textsc{M. Yor}: \textit{Les in\'{e}galit\'{e}s de sous-martingales comme
cons\'{e}quence de la relation de domination}, Stochastics
\textbf{3}, (1979), no.1, 1-15.



\end{thebibliography}
\end{document}